\documentclass[10pt,twocolumn]{article}
\usepackage{geometry}
\usepackage{graphicx}
\usepackage{subfigure}
\usepackage{amsmath, amsthm, amssymb}
\usepackage{fancyhdr}
\usepackage{upgreek,bm,mathtools,amsmath, amsthm, amssymb,graphicx}
\usepackage{graphicx}
\usepackage{setspace,color}
\usepackage{upgreek,bm,mathtools,amsmath, amsthm, amssymb,graphicx}
\usepackage{graphicx}
\usepackage{multirow}
\usepackage{authblk}

\newcommand{\ffrac}[2]{\ensuremath{\frac{\displaystyle #1}{\displaystyle #2}}}
       
\begin{document}
\begin{minipage}{\textwidth}
	\thispagestyle{empty}
	{\bf \Large Cutting force prediction based on a curved uncut chip thickness model\\}
	{\bf David Hajdu, }
	{\bf Asier Astarloa, and}
	{\bf Zoltan Dombovari}
	{\\ \\
		\begin{center}
		\includegraphics[scale=0.8]{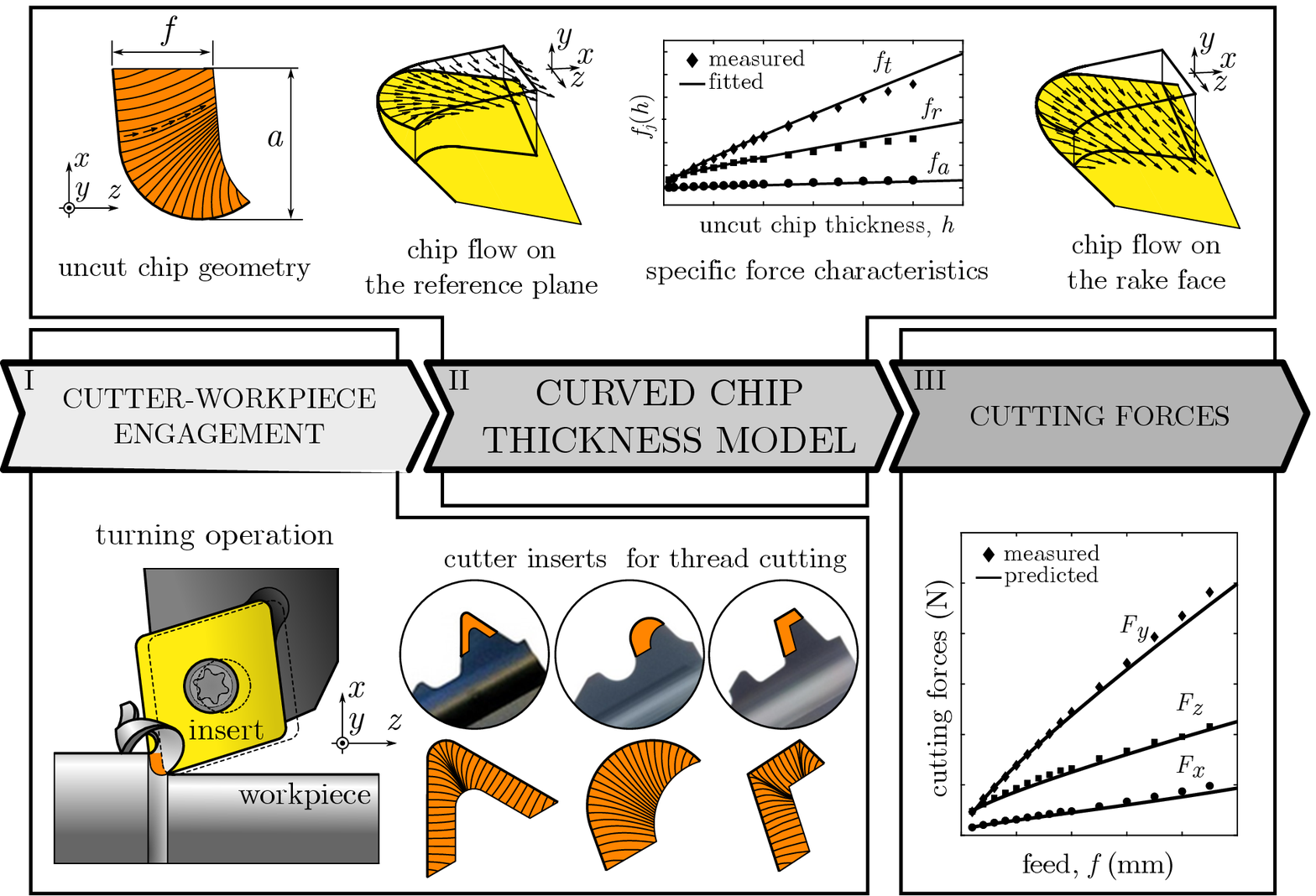}\\
		{graphical abstract}\\
\end{center}}

\end{minipage}
\newpage

\title{\bf Cutting force prediction based on a curved uncut chip thickness model}

\author[1]{\bf David Hajdu\thanks{hajdu@mm.bme.hu}}
\author[2]{\bf Asier Astarloa}
\author[1]{\bf Zoltan Dombovari}

\affil[1]{\small MTA-BME Lend\"{u}let Machine Tool Vibration Research Group, H-1111 Budapest, Hungary\normalsize}
\affil[2]{\small Ideko - Member of BRTA, Dynamics \& Control Department, Elgoibar, Basque Country, Spain\normalsize}

\date{\today}
\maketitle

\setcounter{page}{1}

\begin{abstract}
\it{The curved uncut chip thickness model is presented to predict the cutting forces for general uncut chip geometries. 
{\color{black} The cutting force is assumed to be distributed along a curved path on the rake face of the cutting tool, which makes the solution computable for inserts with nose radius and more complex cutting edge geometries.
	The curved paths originate from a basic mechanical model (a compressed plate model), which is used to mimic the motion of the chip on the rake face of the tool without performing real cutting simulations. 
	Consequently, actual cutting forces are predicted using orthogonal cutting data and the orthogonal-to-oblique transformations.
	The solution satisfies the classical observations and assumptions made on the chip formation process, it is mathematically unique, free of inconsistency and computationally effective.
	Case studies are presented on real cutting tests. 
	The results highlight that cutting force components can be sensitive to modeling assumptions in case of extreme machining parameters.}}
\end{abstract}




\section{Introduction}\label{sec:intro}
Prediction of cutting forces in machining operations is an important part of process planning, since it determines the power consumption of the machine, the required torque during cutting, the static deformation of the tool, surface quality of the workpiece, accuracy, and many other factors.
Accurate modeling of the real mechanics behind the cutting process is, however, very difficult due to the intricate nature of chip formation and physical behavior of the material under cutting conditions. 
Still nowadays, numerical simulations are typically too complicated and time-consuming, or just not reliable enough to replace experimental cutting operations performed on the selected material with the specific tools (see, e.g., finite element models of chip formation processes \cite{Arrazola2010a,Ozel2011a}). 
Therefore, researchers and manufacturing engineers rather rely on the orthogonal cutting database \cite{Altintas2012book}, and predict the cutting forces for general machining operations considering the technological parameters, and the assumptions made on the chip geometry and on the chip flow direction.
These geometric assumptions are, however, inconsistent in general, and it is not trivial which model predicts the real cutting forces better. 

Tools used in practice often have considerable nose radius, since it provides a stronger cutting edge, improves surface quality, and reduces the cutting temperature by the increased heat transfer between the tool and the workpiece \cite{Arsecularatne1995a}. 
It is a favorable property from the technological point of view, but the nose deviates the geometry from the ideal one and makes predictions more complicated and less accurate.
In the current study, we mostly focus on the cutting force prediction for regular single point cutting tools with nose radius, however, the aim of the research is to generalize cutting force predictions for {\color{black} various cutting edge geometries and machining applications.}
In this manner, by using the proposed methodology, the complex geometry of multi-tooth threading inserts{\color{black}, serrated milling cutters and broaching tools} can be dealt generally.

In this paper, a curved {\color{black} uncut} chip thickness model is introduced, which overtakes the most widely used geometrical assumptions.
{\color{black} The solution of the model is based on a basic mechanical model (a compressed plate model), which mimics the motion of the chip on the rake face, without performing real cutting simulations.}
The proposed model can cooperate with orthogonal cutting data and orthogonal-to-oblique transformation, which makes it implementable by using already known empirical cutting models \cite{Altintas2012book}.

The structure of the paper is as follows.
In the rest of Sec.~\ref{sec:intro}, a brief literature review is given on the most commonly used cutting force models, and the orthogonal-to-oblique transformation is presented. Section~\ref{sec:curved_chip} introduces the new approach based on the curved {\color{black}uncut} chip thickness model. 
The geometrical model and comparisons are presented in short case studies. 
Laboratory and industrial experiments are presented in Sec.~\ref{sec:experiment} and Sec.~\ref{sec:industrial}, and finally, results are concluded in Sec.~\ref{sec:conclusion}.

\subsection{Review of existing works}\label{sec:cutforce_literature}

Most of the analytical models dealing with chip formation process in orthogonal cutting were introduced before the 1970s \cite{Germain2013a}. 
These empirical approaches aim to predict the cutting forces based on the essential data only (such as the shear plane angle, friction angle, shear stress, and geometrical data), and omit the complicated phenomena in order to keep the solution tractable. 
The first important model has been presented by Merchant \cite{Merchant1945a,Merchant1945b} for orthogonal cutting  in 1945, where the primary shear zone is a thin plane, and  the cutting edge is perfectly sharp. 
This model allows one to approximate cutting characteristics based on simple algebraic calculations. 
Later, in 1951, Lee and Shaffer \cite{Lee1951a} introduced the slip-line model considering the laws of plasticity, which assumes that the material immediately
above the shear plane is quasi-plastic with slip-lines
being parallel to and perpendicular to the shear plane \cite{Palmer1959a}. 
These models were extended further to describe the real deformation better, e.g., considering the effect of work-hardening, edge radius, strain rate, or heat generation \cite{Jin2011a}. For a more detailed review on the topic, see the work of Germain et al. \cite{Germain2013a}.

The analytical models are not capable to handle the effect of complicated geometries, nose radius or chip breaking grooves \cite{Altintas2012book}.
A mechanistic approach offers an alternative solution, which requires a parameter identification based on cutting experiments performed with the selected workpiece-cutter combination. 
As a result of orthogonal cutting experiments, the cutting force characteristics (as the function of the uncut chip thickness) can be determined.
This methodology makes it possible to predict cutting forces for more complicated tool edge geometries, such as tools with considerable nose radius, {\color{black} however, it does not take into account the effect of chip groves}.

The first model, which considers the effect of nose radius, originates from Colwell \cite{Colwell1954a}, who investigated orthogonal cutting operations already in 1954. 
In that model, it is assumed that the chip flow direction is perpendicular to the equivalent chord connecting the surface and the side points of the edge engagement, while the cutting force magnitude is proportional to the total uncut chip area. Due to the simplicity and applicability of the approach, this assumption has been adopted by many researchers in the past.  Bus et al. \cite{Bus1970a} compared the predictions to experiments, while  Hu et al. \cite{Hu1986a}  and Wen et al. \cite{Wen2003a} extended the method to oblique cutting operations. The assumptions have been implemented by Eynian et al. \cite{Eynian2009a}, who used the equivalent chord to predict unstable self-excited ({\it chatter}) vibrations. 

Although the prediction of Colwell works reasonably well in practice, it is difficult to generally accept that the direction of the chip flow is determined only by the side-points of the engagement edge curve. Instead, Okushima and Minato \cite{Okushima1959a}  suggested in 1959 to use an averaged chip flow direction based on the nearly orthogonal cutting conditions. Since in classical orthogonal cutting, the chip flow direction is normal to the cutting edge, the average of the elementary chip flow angles along the cutting edge can be used to predict the averaged chip flow direction. 

In order to generalize the geometrical models and find a better description of the cutting force direction, Young et al. \cite{Young1987a} proposed a new method in 1987. 
Young's approach suggests that the elementary chip flow direction along the cutting edge is normal to the edge (as in Okushimas's method), the elementary friction force is normal to the edge and the magnitude is determined by the equivalent uncut chip thickness and elementary uncut chip area. 
The solution is, however, not unique, since the geometric
assumptions are only valid close to the cutting edge.
The total uncut chip area is subdivided into such segments, which are bounded by straight lines and perpendicular to the local cutting edge. 
Each segment is viewed as an elementary orthogonal classical cutting, and the resultant force is the sum of elementary cutting forces. 
This approach has been adopted the most in the literature, and found to be reliable for practical applications, however, it fails in case of general {\color{black} uncut chip geometries}. 

The results of Young have been extended by Wang and Mathew \cite{Wang1995a} and later by Arsecularatne et al. \cite{Arsecularatne1995a,Arsecularatne1998a} to tools with nonzero inclinations angles by introducing equivalent cutting edges to trace back the calculation to classical models.
Regenerative machine tool vibrations are predicted based on these approaches by Totis and Sortino \cite{Totis2014a} and by Kuster and Gygax \cite{Kuster1990a} for internal turning, and by Totis for milling operations \cite{Totis2017a}, too.
The predictions have been deeply studied in the literature, experiments are performed and compared to the theoretical models by Lazoglu et al. \cite{Lazoglu2002a}, Atabey et al. {\cite{Atabey2003a,Atabey2003b}}, Kaymakci et al. \cite{Kaymakci2012a}, and by several other researchers.

\begin{figure*}[!htb]
	\centering
	\includegraphics[scale=1.0]{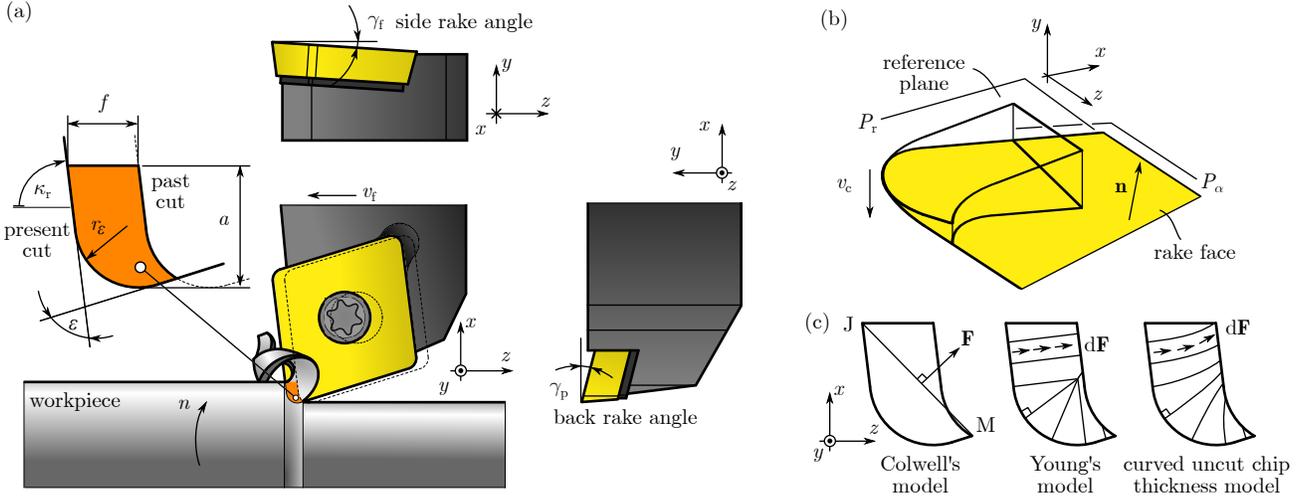}
	\caption{{\bf Cutting force prediction for tools with nose radius in case of turning operation.} (a) Tool geometry. (b) 3D view of the uncut chip area on the reference plane and its projection on the rake face. (c) Direction of the friction force on the reference plane suggested by the different models in case of orthogonal cutting ($\gamma_{\rm p}=0^\circ$ and $\gamma_{\rm f}=0^\circ$).}
	\label{fig:tool}
\end{figure*}

Young's method (and its extensions) also has limitations, e.g., the assumption that the local chip flow direction is always perpendicular to the corresponding cutting edge may produce self-intersecting lines in the model. 
It is not always trivial how to split the area into straight segments (and the physics behind is not ever clear). 
Khoshdarregi and Altintas \cite{Khoshdarregi2015a} presented a method, which keeps the segments straight, but lets the segments tilt to avoid the inconsistent scenarios. 
The method is presented for threading, where the application of Young's method often runs into such difficulties.
The extensions of Young's model, however, still do not guarantee that the solution is mathematically unique.

The differences between the predictions of the above mentioned models have been studied by Wang et al. \cite{Wang2011a}, Kouadri et al. \cite{Kouadri2017a} and Wu et al. \cite{Wu2020a}, just to mention a few. It is often demonstrated, that the calculations are typically acceptable, but measurements are often found to be somewhere in between the suggestions of Colwell and Young. 

In order to approach the solution from a different perspective, nonlinear finite element simulations have also been investigated, see the work of Arrazola and \"{O}zel \cite{Arrazola2010a}, \"{O}zel et al. \cite{Ozel2011a},
Afazov et al. \cite{Afazov2010a} or Wu et al. \cite{Wu2020a}.
However, the complexity of finite element modeling is high, and it is far to fully replace experimental cutting tests, which are often needed to tune the input data in numerical simulations. 

By introducing the curved {\color{black}uncut} chip thickness model, we can cover a wider spectrum of applications and can overcome the difficulties generated by the complicated geometries. 
It is assumed that the cutting force is normal to the cutting edge close to the  edge, but then it bends to avoid self-intersecting directions, mimicking the real chip flow path on the rake face. 
Since the flow of the material is a complicated physical phenomenon, it is not possible to predict accurately the cutting force based on geometrical models only.
Therefore, we create an artificial vector field from a mechanical model, which fulfills the above mentioned criteria and can be an acceptable approximation for the local flow directions. 
The details of the model are presented in Sec.~\ref{sec:curved_chip}.

\subsection{Cutting tool geometry}

An external turning tool with a single point insert is presented in Fig.~\ref{fig:tool}a,b, which is used for demonstration purposes. 
The cutting force models are applicable for milling, drilling, broaching, planing and other non-abrasive machining operations, by modifying the description of geometry and (or) coordinate system.  
The tool moves in the direction of the feed ($z$) with constant speed $v_{\rm f}$, while the feed in one revolution is $f$. The workpiece rotation $n$ is given in rpm, which gives the cutting speed $v_{\rm c}=2\uppi Rn/60$, where $R$ is the radius of the workpiece. The depth of cut $a$ is measured in direction $x$. The main cutting edge angle is $\bf \kappa_{\rm r}$, while $\gamma_{\rm f}$ and $\gamma_{\rm p}$ are the side- and back rake angles of the tool, respectively. The insert has a nose radius $r_{\varepsilon}$, and nose angle $\varepsilon$. 
The coordinate system ($x,y,z$) in Fig.~\ref{fig:tool} is standard for turning operations. 
More specifically, the coordinate system of the machine is denoted by ($x_0,y_0,z_0$), however, the subscripts '$0$' are omitted for simplicity.  
In the example above (Fig.~\ref{fig:tool}), the normal vector of the rake face is determined by the side and back rake angles as
\begin{equation}
	{\bf n}=\begin{bmatrix}
		0\\-\sin\gamma_{\rm f}\\ \cos\gamma_{\rm f}
	\end{bmatrix}\times\begin{bmatrix}
		\cos\gamma_{\rm p}\\-\sin\gamma_{\rm p}\\0
	\end{bmatrix}.
\end{equation}
Considering the geometrical data of the tool (edge angles, nose radius, normal vector, etc.), the edge engagement can be defined. 
Due to simplicity, the chip breakers on the  rake face are omitted, and the surface is assumed to be completely flat (the normal vector is constant along the rake face).

\subsection{Specific cutting force models}\label{sec:o2o}
\begin{figure*}[!htb]
	\centering
	\includegraphics[scale=1]{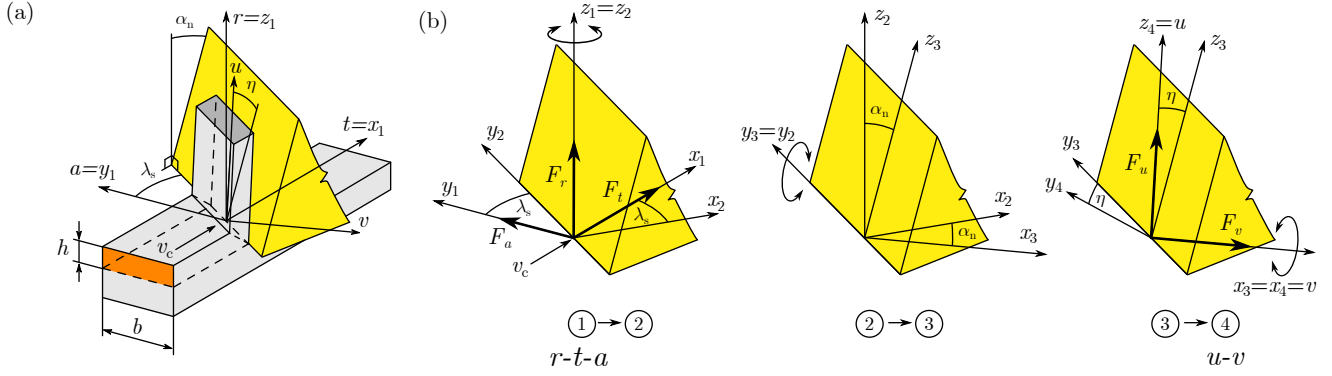}
	\caption{{\bf The orthogonal-to-oblique transformation \cite{Altintas2012book}.} (a) Coordinate systems in oblique cutting operation. (b) Transformation from coordinate system 1 to 4.}
	\label{fig:OrthOb}
\end{figure*}

As it is mentioned above, the new model is based on classical assumptions made on orthogonal and oblique cutting operations. In this short subsection we recall the classical cutting force models and orthogonal-to-oblique transformation (Fig.~\ref{fig:OrthOb}) in order to apply it later for the curved chip thickness model. 
{\color{black} 
	The cutting force components are defined in a coordinate system, which can be oriented as the machine tool, or can be oriented following the direction of the cutting velocity and feed directions. 
	For example, in orthogonal cutting, the tangential ($t$) direction is parallel to the direction of the cutting velocity, while the radial ($r$) and axial ($a$) directions are perpendicular to it.
	However, this might be confusing if different operations are studied at the same time.
	In order to describe the force components concisely, the following coordinate systems (c.s.) and notations  are applied: 
	\begin{itemize}
		\item (0) machine c.s. ($x_0$, $y_0$, $z_0$) or ($x$, $y$, $z$) ; 
		\item (1) orthogonal c.s. $r$-$t$-$a$ ($x_1=t$, $y_1=a$, $z_1=r$); 
		\item (2) cutting edge c.s. ($x_2$, $y_2$, $z_2$); 
		\item (3) rake face c.s. ($x_3$, $y_3$, $z_3$); 
		\item (4) chip flow c.s. $u$-$v$ ($x_4=v$, $y_4$, $z_4=u$).
	\end{itemize}
	The coordinate systems are presented in Fig.~\ref{fig:OrthOb}.}
\noindent Note that the chip flow direction $u$ is defined on the rake face, and $v$ is normal to it. 

Transformation between the different coordinate systems in oblique cutting can be described by the rotation matrices
\begin{equation}
	\begin{split}
		{\bf R}_{z_1,{\lambda_{\rm s}}}  & =\begin{bmatrix}
			\cos\lambda_{\rm s} & -\sin\lambda_{\rm s} & 0\\ \sin\lambda_{\rm s} &\cos\lambda_{\rm s} &0\\
			0 & 0 &1
		\end{bmatrix}, \\
		{\bf R}_{y_2,{\alpha_{\rm n}}}  & =\begin{bmatrix}
			\cos\alpha_{\rm n} & 0 & \sin\alpha_{\rm n} \\ 0 & 1 & 0\\ -\sin\alpha_{\rm n} & 0 &\cos\alpha_{\rm n}
		\end{bmatrix},
		\\
		{\bf R}_{x_3,{\eta}} & =\begin{bmatrix}
			1 & 0 & 0 \\ 0 & \cos\eta & -\sin\eta\\0 & \sin\eta &\cos\eta
		\end{bmatrix},
	\end{split}
\end{equation}
where $\lambda_{\rm s}$ is the inclination angle, $\alpha_{\rm n}$ is the rake angle, and $\eta$ is the chip flow angle measured on the rake face, moreover, the positive direction of rotation is defined around the positive direction of the corresponding axis.
For instance, the transformation between c.s. 1 and 4 is written as
\begin{align}\label{eq:T41_ref}
	{\bf T}_{4,1}={{\bf R}_{{z_1},{\lambda_s}}}{{\bf R}_{{y_2},{\alpha_{\rm n}}}}{{\bf R}_{{x_3},{\eta}}},
\end{align}
where the sequence of multiplications follows intrinsic rotations, i.e., rotations occur in the local frame of reference in the order $z_1$-$y_2$-$x_3$. The indices $(4,1)$ in ${\bf T}_{4,1}$ indicate that ${\bf T}_{4,1}$ transforms the basis vectors of c.s. 1 to the basis vectors of c.s. 4. (see Fig.~\ref{fig:OrthOb}b).

Cutting force components are typically characterized by empirical formulas, which are (linear or nonlinear) functions of at least the uncut chip thickness $h$. 
The specific cutting forces are formulated generally as
${\bf f}(h)={\rm col}f_j(h)$,
where $j=t,a,r$ or $j=u,v$ depending on the coordinate system \cite{Altintas2012book}. 
These characteristics are either determined by orthogonal cutting tests having $f_j(h)$ using the mechanistic approach, or they can be predicted theoretically using material properties and the engineering database \cite{Altintas2012book}.
One of the general model assumes $f_j(h)$ in the form
\begin{equation}\label{eq:cutchar}
	f_{j}(h)= K_{j\rm c}h+K_{j  \rm e}, \;\; j=u,v, \;\;\text{or} \;\; j=t,a,r,
\end{equation}
where $K_{j  \rm c}$ is cutting coefficient and $K_{j  \rm e}$ is the edge coefficient. 
Due to the different nature of the shearing and ploughing (also called ``rubbing'' \cite{Altintas2012book}), the edge force is modeled as a distributed force along the cutting edge, while the friction force and normal force are distributed over the contact area of the rake face. 
Note that the cutting characteristics may appear as a linear function of the uncut chip thickness $h$, however, $K_{j\rm c}(h)$ can also be nonlinear (e.g., titanium-based alloys \cite{Altintas2012book}).
{\color{black} It is an important remark here that the effect of edge radius is not investigated in details, but is assumed be included in the characteristics $f_j(h)$ either by the edge coefficient or as a nonlinear function of the uncut chip thickness.}

The orthogonal-to-oblique transformation is utilized to predict the elementary cutting coefficients $K_{j\rm c}$ in the chip flow c.s. ($j=u,v$) \cite{Altintas2012book,Kaymakci2012a} in the form
\begin{equation}
	\begin{split}
		K_{u \rm c}=\frac{\tau_{\rm s}\sin{\beta_{\rm a}}\sqrt{1+\tan^2\eta\sin^2\beta_{\rm n}}}{\cos\lambda_{\rm s}\sin\phi_{\rm n}\sqrt{\cos^2(\phi_{\rm n}+\beta_{\rm n}-\alpha_{\rm n})+\tan^2\eta\sin^2\beta_{\rm n}}},\\
		K_{v \rm c}=\frac{\tau_{\rm s}\cos{\beta_{\rm a}}\sqrt{1+\tan^2\eta\sin^2\beta_{\rm n}}}{\cos\lambda_{\rm s}\sin\phi_{\rm n}\sqrt{\cos^2(\phi_{\rm n}+\beta_{\rm n}-\alpha_{\rm n})+\tan^2\eta\sin^2\beta_{\rm n}}},
	\end{split}
\end{equation}
where $\rm \tau_{\rm s}$ is the shear stress of the workpiece material, $\phi_{\rm n}$ is the shear angle, $\tan\beta_{\rm n}=\tan{\beta_{\rm a}}\cos\eta$, moreover $\beta_{\rm n}$ is called the projected friction angle, and $\beta_{\rm a}$ is the average friction angle. For simplicity, we follow Stabler's rule, i.e., $\eta=\lambda_{\rm s}$ \cite{Stabler1951a}.


The pressure distribution along the rake face is often modeled empirically, which can be taken into account by a theoretical weight function $\Lambda_j(r,h)$, which satisfies the condition
\begin{equation}\label{eq:Lambda_def}
	\int_{0}^{h}\Lambda_j(r,h) \,{\rm d}r=1, \;\; j=u,v.
\end{equation}
Standard formulas presented in this section express the cutting characteristics in terms of the uncut chip thickness $h$ and the uncut chip width $b$ measured on the reference plane (see Fig.~\ref{fig:OrthOb}a).
Therefore, the cutting force vector ${\bf F}={\rm col} F_j$, $(j=t,a,r)$ is defined theoretically in the form
\begin{equation}\label{eq:cuttingforce_total_int0}
	{\bf F}= \int_{0}^{b} \int_{0}^{h} {\bf T}_{4,0}\begin{bmatrix}
		\Lambda_v(r,h)f_v(h) \\ 0\\ \Lambda_u(r,h)f_u(h)
	\end{bmatrix}{\rm d}r\,{\rm d}a.
\end{equation} 
Since the elementary cutting force vectors are aligned parallel to each other, the integral simplifies to ${ F_j}=bf_j(h)$, where $j=t,a,r$.
Also note that the distribution does not affect the force magnitudes, if the elementary chip segments are straight. 
However, in case of the curved chip model, the weight functions $\Lambda_j$ can have effect. 

%
%

\section{The curved uncut chip thickness model}\label{sec:curved_chip}
The first step to determine the curved chip segments is to generate a vector field, along which the chip is supposed to flow on the rake face. 
This curved chip thickness (see Fig.~\ref{fig:FE_chip}) allows using empirical cutting force characteristics described in Sec.~\ref{sec:o2o}.
For this purpose, we create a hypothetical finite element plate model defined over the uncut chip area with uniform thickness and compress it normal to its plane (in direction $t$ in the $r$-$t$-$a$ c.s.). 
{\color{black} This actually mimics the tangential inertial pressure and its effect on the arriving chip segments.}
The nodal solution is used to generate a vector field for deriving the chip flow path directions and consequently  the geometry of the curved uncut chip thickness.

{\color{black} The consequence of the model is that its solution is
	\begin{itemize}
		\item mathematically unique, because the calculation is based on a linear finite element model, which has only one solution, and
		\item free of inconsistency, i.e., the curved chip segments do not intersect each other.
\end{itemize}}

\subsection{Modeling approach}

\begin{figure*}
	\centering
	\includegraphics[scale=1]{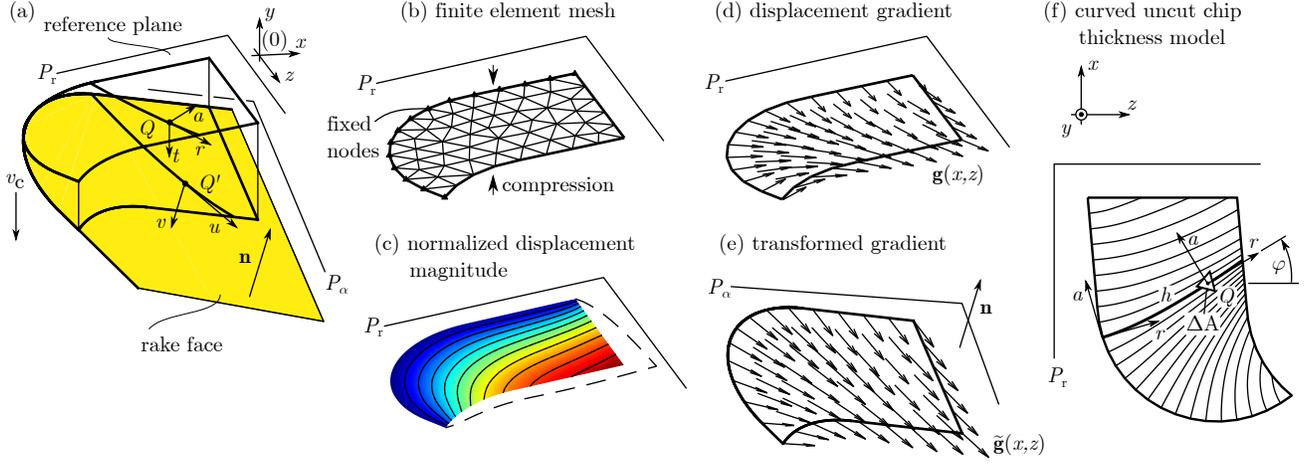}
	\caption{{\bf Generation of the curved uncut chip segments.} (a) 3D view of the uncut chip geometry. (b) Finite element mesh (nodes are fixed on the cutting edge, the load is a compression in direction $y$). (c) Nodal displacement magnitudes calculated from the deformed shape. (d) Displacement gradient ${\bf g}(x,z)$ (gradient is normal to the cutting edge on the reference plane). (e) vector field transformed to rake face following Stabler's rule and the local rake and inclination angles. (f) Curved chip segments generated by the trajectories of the vector field ${\bf g}(x,z)$ on the reference plane.}
	\label{fig:FE_chip}
\end{figure*}

The vector field generating the curved chip thickness on the reference plane is calculated from a linear low-resolution finite element model {\color{black} without any real chip formation simulation}.
We use a simplified isoparametric triangular element, which models a plate with uniform thickness, and only uniform compression normal to its plane is allowed.
The advantage of such a simplified element is the faster numerical computation and easier implementation. 
The displacements are only constrained at the nodes along the cutting edge, while the load is a displacement-driven compression in the normal direction of the reference plane.
The numerical calculation procedure is presented in details in Appendix~\ref{app:FE}, while a sample solution is shown in Fig.~\ref{fig:FE_chip}a-d for turning operations. 
The deformation gradient on the reference plane ($x,z$ in the machine coordinate system) produces the vector field ${\bf g}(x,z)=[g_{x},0,g_z]^\top$, along which the curved uncut chip thickness is defined. 
This approach results in the candidate solution to the curved uncut chip geometry on the reference plane.

If the side and back rake angles are nonzero, then the equivalent rake angles $\alpha_{\rm n}$ and inclination angles $\lambda_{\rm s}$ are different at every point along the rake face. 
Therefore, rotation matrices are applied to describe the transformation from the reference plane to the rake face.
The aim is to treat the cutting process locally as an equivalent oblique cutting.
First, we determine the angle between the gradient vector ${\bf g}$ and axis $z$, denoted by $\varphi$, which is
\begin{equation}
	\tan\varphi=\frac{g_x}{g_z}.
\end{equation}
The corresponding local rake angle and inclination angle are determined from the constraint that the normal vector ${\bf n}=[n_x,n_y,n_z]^\top$ of the rake face is known in the machine c.s..
In order to describe the rotations correctly, the transformation matrix ${\bf T}_{4,0}$ must be determined, which transforms the unit base vector $[1,0,0]^\top$ from the c.s.~0 to 4. This transformed vector is parallel with the normal vector of the plane, i.e.,
\begin{equation}\label{eq:transfangles}
	{\bf T}_{4,0}\begin{bmatrix}
		1\\0\\0
	\end{bmatrix}=\begin{bmatrix}
		-\sin\alpha_{\rm n}\sin\varphi+\cos\alpha_{\rm n}\cos\varphi\sin\lambda_{\rm s}\\
		-\cos\alpha_{\rm n}\cos\lambda_{\rm s}\\
		-\sin\alpha_{\rm n}\cos\varphi-\cos\alpha_{\rm n}\sin\varphi\sin\lambda_{\rm s}
	\end{bmatrix}\equiv-{\bf n},
\end{equation}
where 
\begin{align}\label{eq:T40def}
	{\bf T}_{4,0}={\bf R}_{y_0,\varphi}{\bf R}_{z_1,-\uppi/2}{\bf R}_{z_1,\lambda_{\rm s}}{\bf R}_{y_2,\alpha_{\rm n}}{\bf R}_{x_3,\eta}
\end{align}
is the transformation considering the definition of coordinate systems ({\color{black} see Fig.~\ref{fig:OrthOb}}).
Since the transformation from one normal vector to the other contains the unknown angles, then eq.~\eqref{eq:transfangles} can be solved, that is, $\alpha_{\rm n}$ and $\lambda_{\rm s}$ are calculated as
\begin{equation}\label{eq:transfanglessol}
	\begin{split}
		\sin\alpha_{\rm n} & =n_x\sin\varphi+n_z\cos\varphi,\\
		\sin\lambda_{\rm s} & =\frac{n_z-\cos\varphi\sin\alpha_{\rm n}}{\cos\alpha_{\rm n}\sin\varphi}.
	\end{split}
\end{equation}
Since the angles for relevant technological parameters are relatively small, the $\arcsin$ function gives the correct angles of rotations. 
Note that the solution can be singular if $\sin\varphi=0$ ($\cos\varphi=\pm1$), in this case eq.~\eqref{eq:transfangles} is solved again and the solution \eqref{eq:transfanglessol} is replaced by
\begin{equation}
	\sin\alpha_{\rm n}=n_z\cos\varphi \;\; \text{and} \;\;
	\sin\lambda_{\rm s}  =\frac{-n_x}{\cos\alpha_{\rm n}\cos\varphi}.
\end{equation} 
Also note that if Stabler's rule is followed \cite{Stabler1951a}, then the chip flow angle $\eta$ is equal to the local inclination angle $\lambda_{\rm s}$, which can be substituted into eq.~\eqref{eq:T40def}.
{\color{black} The gradient vector $\bf g$ is then transformed to $\tilde{\bf g}$ following the above mentioned rule (transformation of unit vectors from c.s. 1 to c.s. 4 by ${\bf T}_{4,1}$), or alternatively it can be expressed by the transformation matrix ${\bf T}_{4,0}$ as
	\begin{equation}\label{eq:g-to-gtilde}
		\tilde{\bf g}(x,z):={\bf T}_{4,1}(x,z){\bf g}(x,z)={\bf T}_{4,0}(x,z)\begin{bmatrix}
			0\\0\\1
		\end{bmatrix},
	\end{equation}
	Note, that for orthogonal cutting ($\alpha_{\rm n}=\lambda_{\rm s}=\eta=0$), the transformations simplify to ${\bf T}_{4,0}={\bf R}_{y_0,\varphi}{\bf R}_{z_1,-\uppi/2}$ (and ${\bf T}_{4,1}=\bf I$), and \eqref{eq:g-to-gtilde} degrades to $\tilde{\bf g}={\bf g}$, as expected.}

The {\color{black}uncut} chip thickness is equal to the length of the streamline measured on the reference plane, and can be calculated numerically by solving the differential
\begin{equation}
	\frac{{\rm d} x}{{\rm d} z}=\frac{\tilde{g}_x(x,z)}{\tilde{g}_z(x,z)}.
\end{equation}
Many numerical softwares provide built-in algorithms to generate streamlines from arbitrary starting points if the vector field is known (e.g., the \texttt{streamline} or \texttt{stream2} functions in MATLAB). The vector field ${\bf g}(x,z)$ must be interpolated between the node points, then numerical solutions can be initiated from the centers of the triangular elements. 
Finding the starting points attached to the cutting edge and end points located on the free boundaries, the length of the path can be determined, which is identical to $h$.


Since the cutting force is distributed along the curved chip thickness, the total cutting force is obtained by summing up the elementary forces over the entire uncut chip area $A$. 
In this manner, the cutting force is theoretically written as
\begin{equation}\label{eq:cuttingforce_total_int}
	{\bf F} =  \int_{(A)} {\bf T}_{4,0}\begin{bmatrix}
		\Lambda_v(r,h)f_{v}(h) \\ 0\\ \Lambda_u(r,h)f_{u}(h)
	\end{bmatrix} {\rm d}A,
\end{equation} 
where ${\bf T}_{4,0}$ is the local transformation matrix, and $h$ is the {\color{black}uncut} chip thickness of the elementary curved chip passing through the infinitesimal chip area ${\rm d}A$.
Note, that the edge coefficients are often modeled as concentrated forces along the edge only. 
Therefore, by losing generality, we assume that the weight function is of the form
\begin{align}\label{eq:Lambda_def2}
	\Lambda_j(r,h):=\frac{K_{j\rm c}(h)+K_{j\rm e}\delta(r-0)}{K_{j\rm c}(h)h+K_{j\rm e}}, \;\; j=u,v,
\end{align}
{\color{black} where $\delta(r-0)$ is the shifted Dirac-delta function, which indicates that the edge force is defined at the positive side of zero}.
This definition automatically satisfies \eqref{eq:Lambda_def}, moreover
\begin{align}
	\int_{0}^h \
	\Lambda_j(r,h)f_{j}(h)\,{\rm d}r=f_j(h)
\end{align}
also holds.
Inserting definition \eqref{eq:Lambda_def2} into eq.~\eqref{eq:cuttingforce_total_int} yields
\begin{equation}\label{eq:F_example}
	{\bf F} =  \int_{(A)} {\bf T}_{4,0}\begin{bmatrix}
		K_{v\rm c}(h) \\ 0\\ K_{u\rm c}(h)
	\end{bmatrix} {\rm d}A+\int_{(L)} {\bf T}_{4,0}\begin{bmatrix}
		K_{v\rm e} \\ 0\\ K_{u\rm e}
	\end{bmatrix} {\rm d}L,
\end{equation}
where the second integral is only evaluated along the cutting edge of length $L$. Note, that formula \eqref{eq:F_example} is directly obtained, if the the edge coefficient is not considered as part of the characteristics $f_j(h)$.
{\color{black} Since the computation of the cutting force defined in the form \eqref{eq:F_example} is not possible analytically in general cases, the numerical integration is written instead in the form}
\begin{equation}
	{\bf F} = \sum_{i=1}^{N_{\rm e}} {\bf T}_{4,0,i}\begin{bmatrix}
		K_{v,\rm c}(h_i)\\0\\K_{u,\rm c}(h_i)\\
	\end{bmatrix} {{\Delta}A_i}+\sum_{j=1}^{N_{\rm k}}{\bf T}_{4,0,j}\begin{bmatrix}
		K_{v\rm e}\\0\\K_{u\rm e}\\
	\end{bmatrix} {\Delta}L_j,
\end{equation}
where $N_{\rm e}$ is the number of discrete elements on the area, $h_i$ is the length of the uncut curved chip segment passing though the center of element numbered $i$, ${\Delta}A_i$ is the area portion, ${\bf T}_{4,0,i}$ and ${\bf T}_{4,0,j}$ are the local transformation matrices assigned to the area and edge segments, ${\Delta}L_j$ is the length of the elementary cutting edge, and $N_{\rm k}$ is the number of the elementary edges.
The numerical integration is based on the triangulation of the uncut chip area, which is already available due to the finite element meshing (see Appendix~\ref{app:FE}).

\subsection{Theoretical comparisons}\label{sec:example}
{\noindent \it{Example 1 - V-shaped insert}}

\begin{figure}[!htb]
	\centering
	\includegraphics[scale=1]{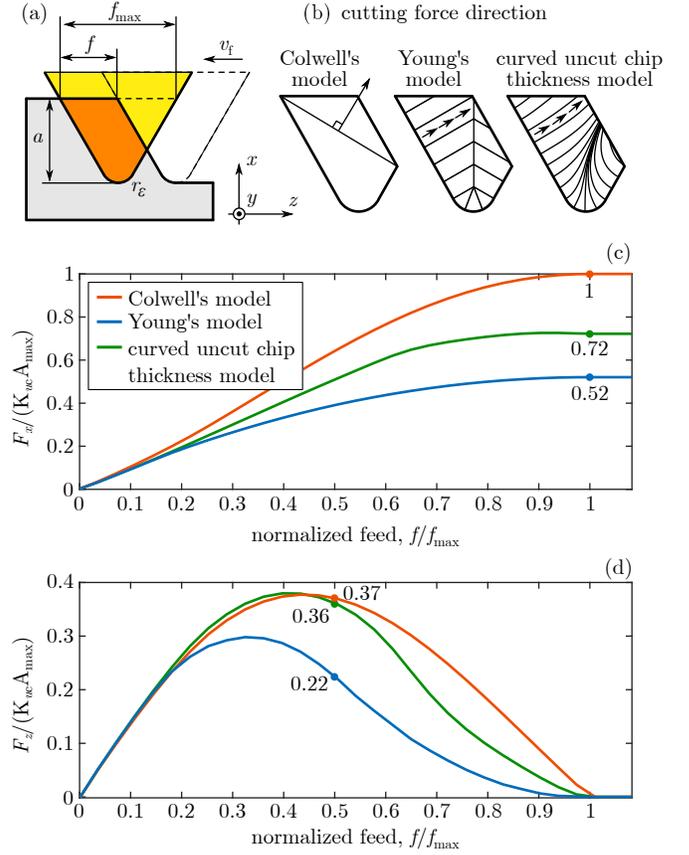}
	\caption{{\bf Theoretical comparison of cutting force models for an orthogonal cutting operation with a V-shaped tool}. (a) Cutting configuration. (b) Cutting force models. (c-d) Normalized cutting forces in direction $x$ and $z$. The maximum discrepancy is approximately 25\%. }
	\label{fig:example}
\end{figure}

{\color{black} The selection of the tool and the holder for turning operations significantly depends on the type of process, material and machine capabilities.
	In order to visualize the predictions of the cutting force models, we present a theoretical example performed with a V-shaped tool for extreme technological parameters that can arise in roughing operations. 
	The assumed data are } $r_{\varepsilon}=0.2$ mm, $\kappa_{\rm r}=60^\circ$, $\varepsilon=60^\circ$, $a=1$ mm, $\gamma_{\rm f}=0^\circ$, $\gamma_{\rm p}=0^\circ$, and the feed goes up to the maximum $f_{\rm max}$, until the subsequent cuts do not overlap each other. 
The geometry of the uncut chip area is seen in Fig.~\ref{fig:example}a.

{\color{black}For simple calculation purposes}, we assume that the cutting forces are linear functions of the uncut chip thickness $h$, and the edge forces are completely omitted, i.e., $f_u(h)=K_{u\rm c}h$ and $f_v(h)=K_{v\rm c}h$. 
{\color{black} Colwell's approach (see \cite{Colwell1954a}) assumes that the magnitude of the cutting force is proportional to the total chip area $A$, while the direction is perpendicular to the equivalent chord.}
This results in the cutting force magnitudes $F_u=K_{u\rm c}A$ and $F_v=K_{v\rm c}A$. 
The maximum force is reached, when the uncut chip area is the largest $A_{\rm max}$. 
In order to reduce the number of parameters, the cutting forces are normalized by $K_{u\rm c}A_{\rm max}$, then the dimensionless feed $f/f_{\rm max}$ is introduced, where $f_{\rm max}$ is the largest feed, which can be reached for certain depth of cut and geometrical parameters (until subsequent cuts do not overlap). The normalized forces are compared accordingly.  
Figure~\ref{fig:example}b presents the elementary chip segments and panels c-d show the normalized cutting forces in direction $x$ (depth of cut) and $z$ (feed). 
The cutting force in direction $y$ for orthogonal cutting (and for the parameters above) is the same for all three models. 
{\color{black} As it can be seen, the cutting model of Young (see \cite{Young1987a}) is difficult to adopt}, because the intersections of the lines cannot be avoided if they are perpendicular to the cutting edge. 
It is not trivial how the model should be tuned to avoid the inconsistent intersections and keep other assumptions. 
As opposed to this, the curved uncut chip thickness model generates curved chip segments, which start perpendicular to the edge and assumingly deflect as the chip flows avoiding the possibility of intersections. This  property is assured by the nature of continuous vector fields.

It can be concluded that the predicted cutting forces are similar at small feeds, but significantly different for large feeds (difference is nearly 25\%). 
In case of internal turning (boring) operations, the deflection of the tool and precision of the machining are mostly affected by the radial force $F_x$.
Therefore, it has to be predicted accurately, since it affects the most the vibratory behavior during roughing and surface quality of finishing operations.
\\

\begin{figure}[!htb]
	\centering
	\includegraphics[scale=1]{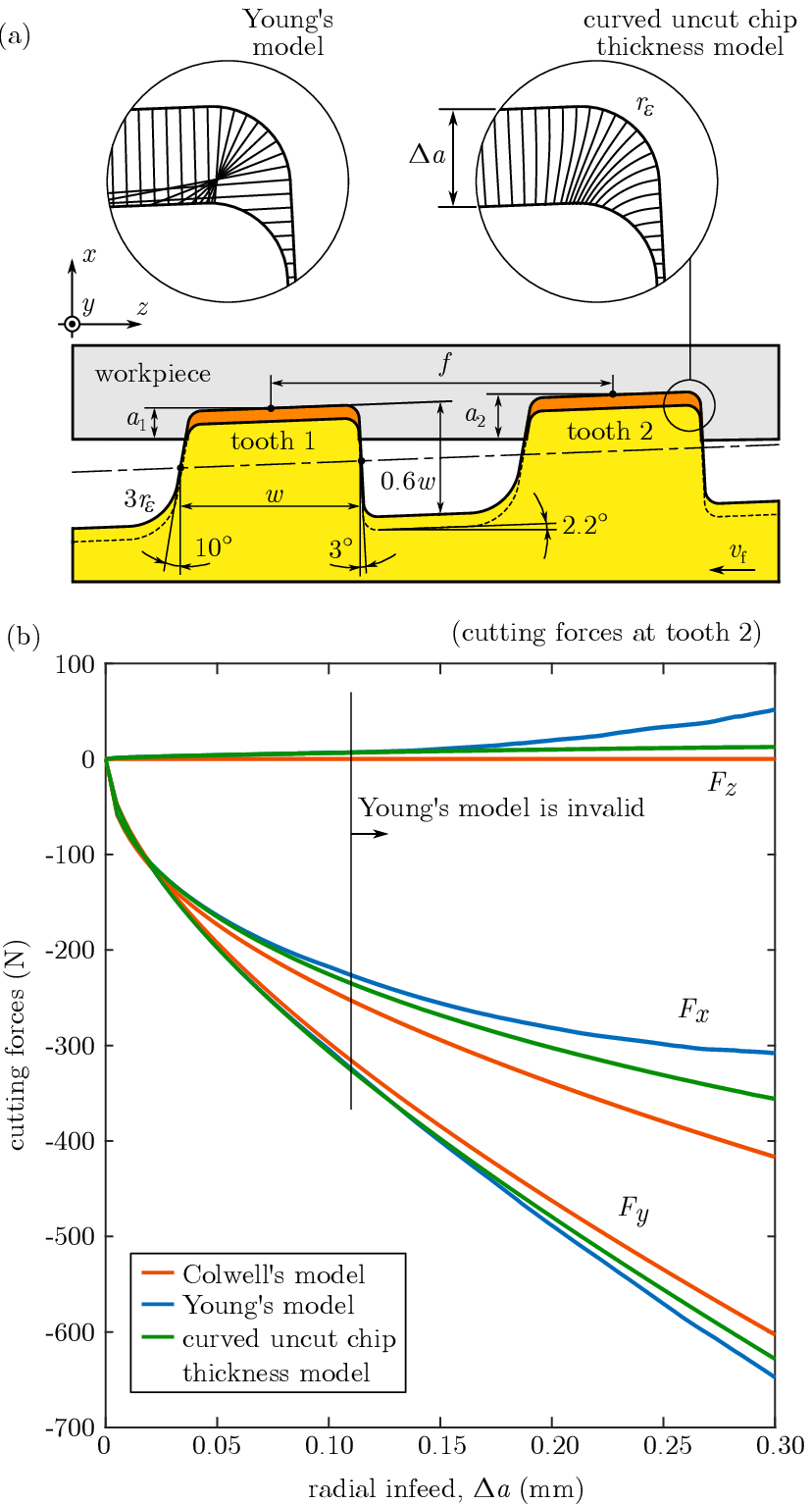}
	\caption{{\bf Theoretical comparison of cutting force models for threading operation}. (a) Two-point threading operation with buttress profile. (b) Predicted cutting forces acting on tooth 2.}
	\label{fig:example2}
\end{figure}

{\noindent \it{Example 2 - threading insert}}

There are other cutting operations, where the uncut chip area formed between subsequent cuts is irregular, such as threading \cite{Khoshdarregi2015a} or gear hobbing \cite{Zheng2021a}.
In such cases, the classical models may not work in all cases, and the original assumption of Young \cite{Young1987a} must be violated, as presented by Khoshdarregi and Altintas \cite{Khoshdarregi2015a} for threading.

A two-point threading insert (buttress) is presented in Fig.~\ref{fig:example2}a, where such scenarios can occur, especially at high radial infeeds $\Delta a$. 
The curved chip segments can be calculated by the presented method easily, which makes the transition to the uncut curved chip thickness smooth at the corners, where discrepancies are expected. 

The following numerical data are assumed during the calculation: $w=1$ mm ($f=2w$ per revolution), $r_\epsilon=0.1$ mm, $a_1=0.1$ mm, $a_2=a_1+\Delta a$.
The work material is assumed to be AISI 1045 steel. The corresponding empirical cutting coefficients are $K_{u \rm c}(h)=691.6h^{-0.534}$ MPa, $K_{v \rm c}(h)=1204.3h^{-0.384}$ MPa, where the numerical data are taken from the work of Khoshdarregi and Altintas \cite{Khoshdarregi2015a} ($h$ is given in mm).
For simplicity, we assume that the cutting operation is orthogonal, the edge coefficients are neglected, the tool is flat, moreover, we only focus on the cutting force acting on tooth 2 (the method can be repeated for the additional teeth).

The case study is presented in Fig.~\ref{fig:example2}, where the depth of cut $a_1=0.1$ mm is fixed, while $\Delta a$ is a design parameter considered up to $0.3$ mm.
The large radial infeeds produce self-intersections according to Young's model \cite{Young1987a} (Fig.~\ref{fig:example2}a), which must be handled with special care during computation. 
In order to avoid unrealistic cases, the straight chip segments were chopped off, if self-intersections are found. 
Note that this simplification makes it possible to calculate the cutting forces, but such solution is mathematically not unique, and larger discrepancies are expected. 
Colwell's model \cite{Colwell1954a} in this case is more unrealistic, since the equivalent chord is parallel with the direction of the feed, providing no force in direction $x$.

The predicted cutting forces are presented in Fig.~\ref{fig:example2}b.
Larger discrepancies between the models are visible at higher infeeds, where the {\color{black}uncut} chip thickness is also larger. 
The nonlinear cutting force characteristics also elevates the differences between the predictions, which support the relevance of the uncut curved chip thickness model.

\section{Laboratory experiments}\label{sec:experiment}

Laboratory experiments are carried out for a titanium-based alloy \mbox{Ti$_6$Al$_4$V} (grade 5), which is a widely used (and fairly expensive) material in aerospace industry.
Machining of these materials is very complicated due to its high chemical activity and poor thermal conductivity resulting high cutting temperatures \cite{Jawahir2016a}. 
In order to extend tool life, tool manufacturers recommend different cutting speed zones for different machining operations, for instance, $v_{\rm c}=$ 20-25 m/min for first stage machining (e.g., heavy roughing or skin removal). However,  in some studies the cutting speed is even below this, only 3 m/min \cite{Aksu2016a}.
Nowadays, for intermediate stage machining operations (e.g., profiling), the cutting speeds are in the range $v_{\rm c}= $ 40-80 m/min, but these operations are often performed with coolant. 
In the experiments presented in the paper, the cutting speed was $v_{\rm c}=21$ m/min for all of the tests (to reach as high feeds as possible) and no cooling was applied. 
{\color{black} A schematic of the cutting test experiments can be seen in Fig.~\ref{fig:Measurement_config}.}

\subsection{Identification of cutting force characteristics}
\begin{figure}
	\centering
	\includegraphics[scale=1]{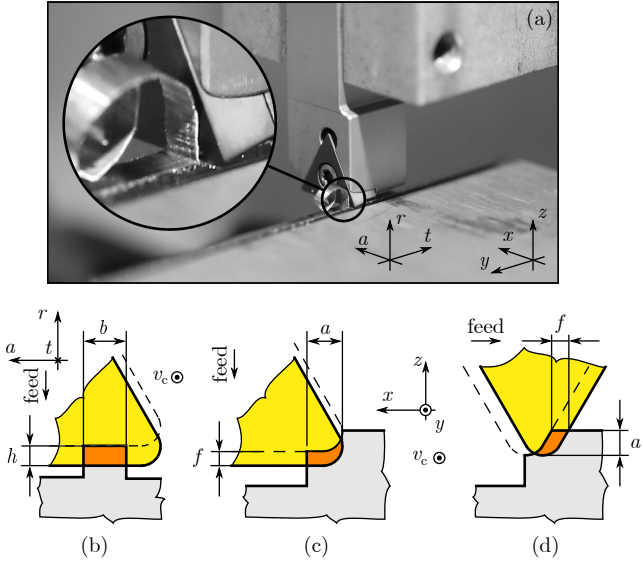}
	\caption{{\bf Schematic of the cutting test experiments.} (a) Chip formation in case of orthogonal cutting on \mbox{Ti$_6$Al$_4$V}. (b) Ringe cutting test with the straight edge of the tool. (c-d) Cutting tests with the nose of the tool.}
	\label{fig:Measurement_config}
\end{figure}

\begin{figure}
	\centering
	\includegraphics[scale=1]{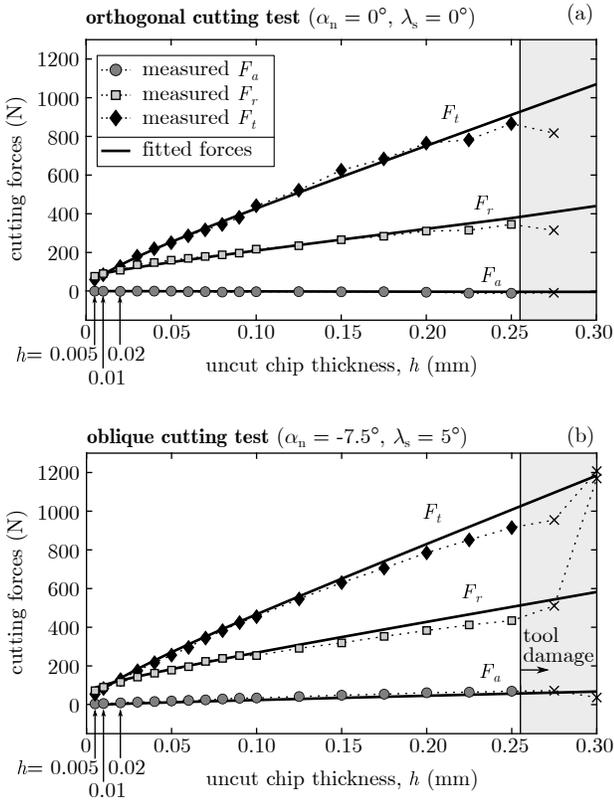}
	\caption{{\bf Orthogonal and oblique ridge cutting tests  on \mbox{Ti$_6$Al$_4$V}}. (a) Measured and fitted force components in orthogonal cutting. (b) Measured and fitted force components in oblique cutting.}
	\label{fig:Measurement1}
\end{figure}

The cutting force characteristics have been identified from cutting tests performed on premachined ridges (see Fig.~\ref{fig:Measurement_config}b). {\color{black} The width of the rigde was $b=2$ mm, while the uncut chip thickness $h$ was varied non-uniformly from $0.005$ to $0.3$ mm (e.g., $h=$ 0.005, 0.01, 0.02, 0.03, 0.04, 0.05, 0.06, 0.07, 0.08, 0.09, 0.1, 0.125, 0.15, 0.175, 0.2, 0.225, 0.25, 0.275 mm).
	The measured cutting forces are plotted in Fig.~\ref{fig:Measurement1}a.}
The data of the insert is given in Table~\ref{table:cut_par}. The cutting forces were measured by a Kistler 9129AA multicomponent dynamometer and data were acquired by NI-9234 input modules and an NI cDAQ-9178 chassis.

\begin{table}
	\caption{{\bf Orthogonal cutting parameters used for titanium-based alloy \mbox{Ti$_6$Al$_4$V}}. $\beta_{\rm a}$ and $\alpha_{\rm n}$ are given in radians \cite{Altintas2012book}.} \label{table:cut_par}
	\begin{tabular}{ll}
		\hline
		Orthogonal cutting data &Insert\\
		\hline
		$\tau_{\rm s}$ = 613 MPa & Widia \\
		$\beta_{\rm a}{(\rm rad)}$ = 0.34 + 0.441$\alpha_{\rm n}(\rm rad)$ & TCMW16T304 THM\\
		$C_0$ = 0.88+0.63$\alpha_{\rm n}(\rm rad)$ & $r_{\varepsilon}= 0.4$ mm\\
		$C_1$ = 0.35+ 0.12$\alpha_{\rm n}(\rm rad)$ & uncoated carbide \\
		$K_{u\rm e}$ = 32.34 N/mm & flat surface (no chipbreaker)\\
		$K_{v\rm e}$ = 4.6 N/mm & \\
		\hline
	\end{tabular}
\end{table}

Each test was repeated four times in order to make sure that the static deformation of the tool is compensated and the actual uncut chip thickness is identical to the prescribed value. 
This way the tool wear and damage could also be monitored.
If the forces change radically, it can be an indicator of tool edge damage.
The uncut chip thickness was also changed randomly between different test cuts below \mbox{$0.2$ mm}. The mean forces show a smooth monotonic tendency indicating that the order of test cuts did not affect the results. Above \mbox{$0.2$ mm} the measured forces were evaluated cautiously, because the heat generated during cutting was large and the edge damage initiated suddenly. 
{\color{black}
	Cutting tests were regularly repeated at $h=0.05$ mm.
	When tool wear initiates, the cutting force components change (sometimes radically) and this change is monitored by repeating the same tests at a fixed uncut chip thickness.
	If the mean forces were approximately the same (below a threshold of $5\%$), then the experiment was continued. 
	When the difference was larger than the limit, the insert was replaced or the test series was completed. 
	This was the limiting uncut chip thickness that could be reached during the measurements. 
	The largest uncut chip thickness, above which the tool wear was always significant are indicated by crosses in Fig.~\ref{fig:Measurement1}a,b. 
}

The same measurement procedure has been conducted for a different tool holder (with the same insert) in order to analyze the variation of cutting parameters. The rake and inclination angles of the setup were $\alpha_{\rm n}=-7.5^\circ$ and $\lambda_{\rm s}=5^\circ$. The cutting forces are shown in Fig.~\ref{fig:Measurement1}b.

The force model presented by Altintas \cite{Altintas2012book} for titanium-based alloys was adopted, and the parameters were tuned according to the measurements. The shear angle $\phi_{\rm n}(h)$ was assumed to depend on the uncut chip thickness $h$. 
The chip compression ratio for this alloy is assumed to be of form $r_{\rm c}=C_0h^{C_1}$, where $h$ is given in mm, and $C_0$ and $C_1$ are empirical parameters depending on the rake angle. 
The shear angle is calculated by the formula
\begin{equation}
	\phi_{\rm n}=\arctan\frac{r_{\rm c}\cos\alpha_{\rm n}}{1-r_{\rm c}\sin{\alpha_{\rm n}}},
\end{equation}
which makes the cutting coefficients $K_{u\rm c}(h)$ and $K_{v\rm c}(h)$ dependent on $h$.
The fitted curves in Fig.~\ref{fig:Measurement1} were calculated using the classical orthogonal-to-oblique transformations presented in Sec.~\ref{sec:o2o}.

\begin{figure*}
	\centering
	\includegraphics[scale=1]{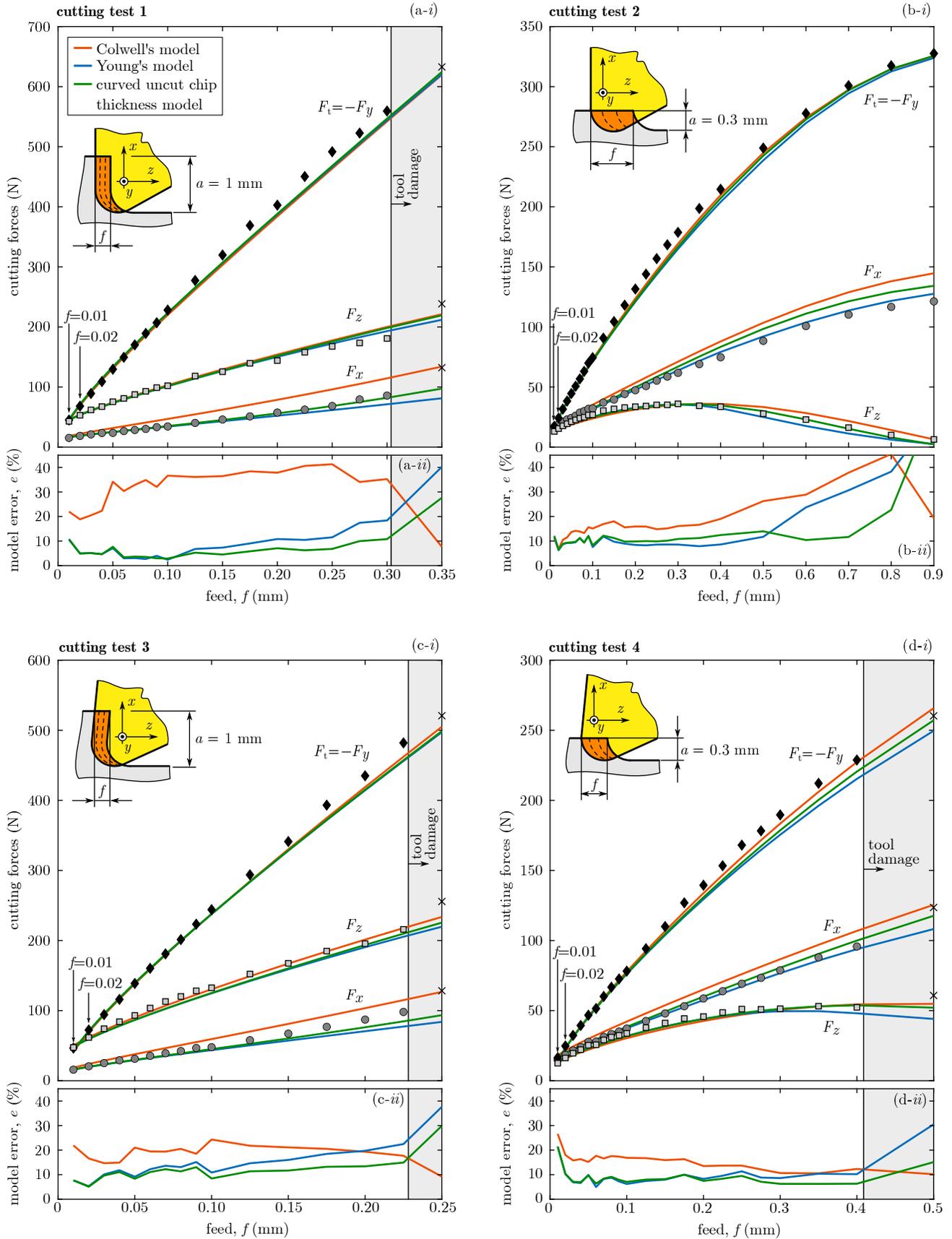}
	\caption{{\bf Cutting force validation on \mbox{Ti$_6$Al$_4$V} with a tool with nose radius (insert type TCMW16T304).} (a-b) Orthogonal cutting tests ($\kappa_{\rm r}=90^\circ$, $\gamma_{\rm f}=0^\circ$, $\gamma_{\rm p}=0^\circ$, $\varepsilon=60^\circ$, $r_{\varepsilon}=0.4$ mm). (c-d) Oblique cutting tests ($\kappa_{\rm r}=93^\circ$. $\gamma_{\rm f}=-7.5^\circ$, $\gamma_{\rm p}=-5^\circ$, $r_{\varepsilon}=0.4$ mm). Colwell's method (red), Young's method (blue) and the curved uncut chip thickness model (green) predict well the cutting force components. The cutting force component $F_x$ is the most sensitive to modeling concepts.}
	\label{fig:Measurement2}
\end{figure*}

\begin{figure*}[!htb]
	\centering
	\includegraphics[scale=1]{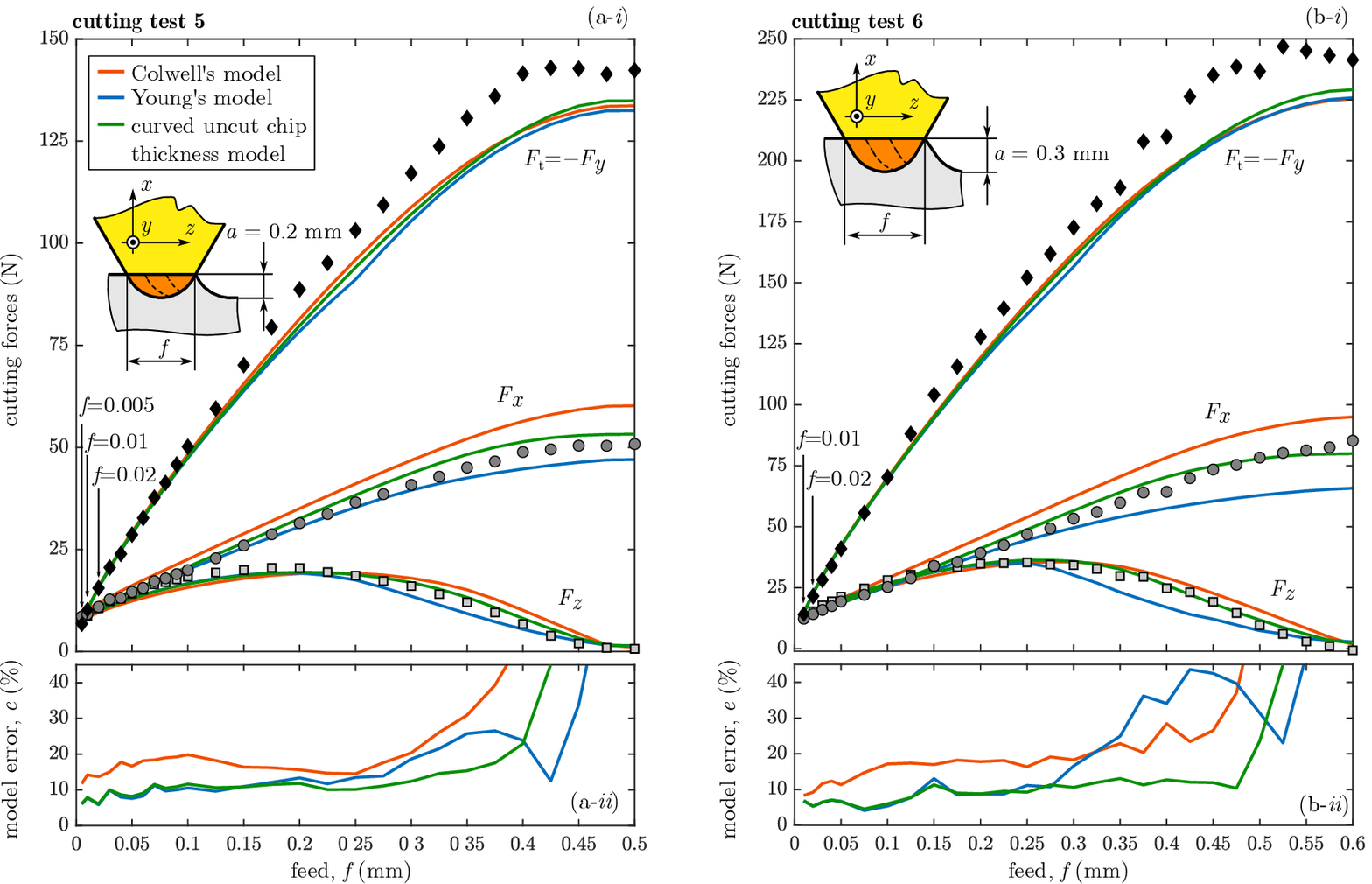}
	\caption{{\bf Cutting force validation on \mbox{Ti$_6$Al$_4$V} with a metric threading insert (16ER1.5ISO).} $\kappa_{\rm r}=60^\circ$, $\gamma_{\rm f}=-0.8^\circ$, $\gamma_{\rm p}=1.2^\circ$, $\varepsilon=60^\circ$, $r_{\varepsilon}=0.216$ mm. The depth of cuts are (a) 0.2 mm and (b) 0.3 mm. The deviations in the cutting fore components $F_x$ and $F_z$ are larger at higher depth of cuts, see Fig.~\ref{fig:example}.}
	\label{fig:Measurement3}
\end{figure*}

\subsection{Cutting tests with nose radius}
Two orthogonal and two oblique cutting test series have been conducted on the same experimental configuration, see Fig. \ref{fig:Measurement_config}c.
The procedures were the same, the tests were stopped when the tool wear initiated. 

The static compliance of the tool holder design was important in these tests, because the deformation in direction $x$ results in change in the actual depth of cut, as opposed to the ridge tests. 
The cuts were repeated four times with the same feed in direction $z$, therefore the deformation in this direction was compensated. The static compliance $\bf C$ of the tool at the tool tip was identified from dynamical impact tests using a modal hammer (Endevco 2302-10), accelerometers (PCB 352C23), and the NI data acquisition system.
The measured compliance matrix of the assembled tool and clamping device corresponding to the tool tip in the ($x,y,z$) coordinate system (see Fig.~\ref{fig:Measurement_config}a) is 
\begin{equation}
	{\bf C}=\begin{bmatrix}
		0.1405   &-0.0172  & -0.0834 \\
		-0.0172    &0.0067   & 0.0076\\
		-0.0834    &0.0076   & 0.0531
	\end{bmatrix}\cdot 10^{-6}\,\frac{\rm m}{\rm N}.
\end{equation}
The actual depth of cut can be calculated by solving a nonlinear equation 
\begin{equation}
	\Delta{\bf r}={\bf C}\,{\bf F}(a-\Delta{ r_x})\\
\end{equation}
using Newton-Rhapson's method, where $\Delta{\bf r}=[\Delta r_x, \Delta r_y,\Delta r_z]^\top$ is the vector of static deformations of the tool tip resulted by the cutting force $\bf F$. Since $\bf F$ depends on the modified depth of cut $\tilde{a}=a-\Delta{ r_x}$, the equation is nonlinear, and must be solved iteratively. 
The stiffness is large, and a few iterations (5-10) converge. The iterations resulted in 5-10\% difference in the cutting force, however, this small change also affects the final characteristics and the comparisons.

The predicted and the measured cutting forces are presented in Fig.~\ref{fig:Measurement2}. Panels a-b show two orthogonal cutting operations, where the prescribed depth of cuts are $a=1$ mm and $a=0.3$ mm, respectively. The feed was increased until the tool wear initiated and the procedure could not be continued. 
The same cutting operation was repeated with an oblique tool holder, the results are shown in panels c-d. In each case, the predicted forces ({\color{black} according to Colwell's model \cite{Colwell1954a}, Young's model \cite{Young1987a}} and the curved uncut chip thickness model) were compared to the measured forces. The tangential force ($F_{\rm t}=-F_{y}$) is the largest, the models predict similar cutting forces, but the measurement tends to give larger forces in magnitude. In each case, the radial force ($F_x$) is the most sensitive to modeling approaches, which has huge impact on the machine's precision and static deformation in case of (internal) turning operations.

In order to visualize and quantify the difference between measurements and models, a normalized modeling error is introduced as
\begin{equation}
	\begin{split}
		{e}:=\sqrt{\sum_{i=x,y,z}\left(\frac{F_{ i\rm, model}-F_{ i\rm, measured}}{F_{ i\rm, measured}}\right)^2},
	\end{split}
\end{equation}
which is the norm of the relative model error. 
The modeling errors are visualized in sub-figures ($ii$) in Fig.~\ref{fig:Measurement2}.
Due to the limited feeds and depth of cuts, it was not possible to reach higher material removal rates and show significant differences without tool damage. 
The last cut at the highest feeds resulted in significant tool wear, which prevented further tests. 
The difference between the models is moderate, the predicted  forces $F_y$ and $F_z$ are typically less distinguishable compared to the measured forces, however, the force component $F_x$ is sensitive (both in relative and absolute value). 
In these experiments Colwell's approach is less accurate, but Young's method and the new curved chip thickness model gives just a minor difference. 

{\color{black} Since the tool damage could not be prevented at higher feeds, it was not possible to continue the experiments on titanium alloy with the cutting insert TCMW16T304. 
	Instead, a 16ER1.5ISO external metric threading insert has been applied, which has a nose radius $r_\varepsilon=0.216$ mm.
	With this tool, the configuration presented in {\it Example 1} (Se.~\ref{sec:example}) can be tested.
	The surface of the insert was flat, however, the normal vector was not complete parallel with the feed direction, meaning that the cutting operation is not purely orthogonal. Since the cutting force components were small, the cutting edge angles were measured, and approximated as $\gamma_{\rm f}=-0.8^\circ$ and $\gamma_{\rm p}=1.2^\circ$.
	The predicted and measured cutting forces are presented in Fig.~\ref{fig:Measurement3}. 
	It must be noted, that the predicted cutting force components $F_x$ and $F_z$ give good agreement with measurements, but the tangential forces tend to underestimate the actual cutting forces.
	In the cutting test series, the depth of cuts $a$ were $0.2$ and $0.3$ mm. At $a=0.4$ mm the tool tip damaged and experiments could not be continued.
	
	\section{Industrial experiments}\label{sec:industrial}
	Industrial experiments have been carried out on a lathe with a turning insert of type DCMW11T304 ($r_\epsilon=0.4$ mm, $\epsilon=55^\circ$).
	The cutting speed was $v_{\rm c}= 300$ m/min, the work material was aluminum 7075 T6, for which the corresponding linear cutting coefficients have been taken from the orthogonal cutting database assuming orthogonal cutting conditions \cite{Dombovari2010a}, i.e., $K_{u\rm c}=229$ MPa and $K_{v\rm c}=856$ MPa.  
	The fitted edge coefficients are $K_{u\rm e}=87$ N/mm and $K_{v\rm e}=15$ N/mm.

	Three sets of experiments (cutting tests 7-9 in Fig.~\ref{fig:Measurement4}) have been carried out, where the radial depth of cuts $a$ were $0.4$, $0.8$ and $1.2$ mm, respectively.
	In each case the feed per revolution $f$ was increased until $f_{\rm max}$ was reached.
	Due to the  extreme cutting parameters, strong chip segmentation appears in most of the cutting tests. 
	The tangential forces ($F_{\rm t}=-F_y$) are the same for all the theoretical models, because the characteristics are linear functions of the uncut chip thickness, the cutting operation is orthogonal, and the edge force is modeled as a distributed force along the cutting edge only.  
	The radial force $F_x$ is still the most sensitive to the modeling concepts, however, the predictions are not accurate enough to tell confidently which model is the most reliable. 
	Strong chip segmentation and inaccurate cutting force parameters can be responsible for the discrepancies.

	\begin{figure}[!htb]
		\centering
		\includegraphics[scale=0.9]{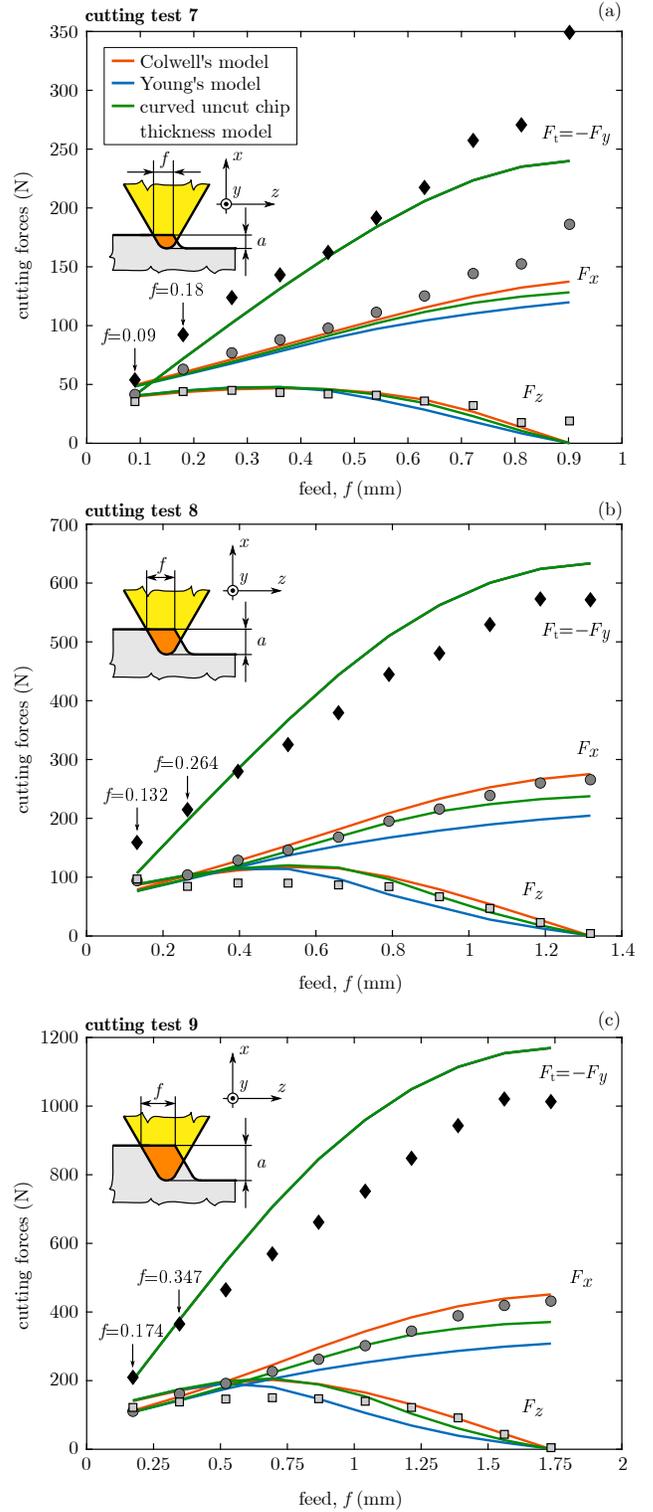}
		\caption{{\bf Industrial validation of predicted cutting forces on AL7075 (insert type DCMW11T304).} ($\kappa_{\rm r}=62.5^\circ$, $\gamma_{\rm f}=0^\circ$, $\gamma_{\rm p}=0^\circ$, $\varepsilon=55^\circ$, $r_{\varepsilon}=0.4$ mm). (a) $a=0.4$ mm, $f_{\rm max}= 0.902$ mm; (b) $a=0.8$ mm, $f_{\rm max}=1.318$  mm; (c) $a=1.2$ mm, $f_{\rm max}=1.735$  mm.}
		\label{fig:Measurement4}
	\end{figure}
}

\section{Conclusion}\label{sec:conclusion}
The curved uncut chip thickness model is presented to generalize the cutting force prediction for arbitrary {\color{black}uncut} chip geometries. The uncut chip area is subdivided into curved segments, along which the elementary uncut chip thickness is defined. 
{\color{black} This semi-analytical model gives an intermediate tool to generate cutting forces for general uncut chip overlaps avoiding any use of large scale finite element solution. The methodology serves a good base for application to cutting operations, where general uncut chip geometries have importance, like milling or broaching.
	The main advantages of the new model can be summarized as follows:
	\begin{itemize}
		\item The curved uncut chip segments are generated from a basic mechanical model (a compressed plate model), which has a unique mathematical solution.
		\item The assumption of Young made on the chip formation process remains valid \cite{Young1987a} (the elementary chip segments are perpendicular to the cutting edge).
		\item The curved uncut chip segments do not overlap each other, the solution is free of inconsistency.
		\item The curved uncut chip thickness model cooperates with the orthogonal-to-oblique transformation and with general empirical cutting force characteristics.
		\item Various cutting edge and uncut chip geometries  can  easily be modeled.
		\item The computation is fast and numerically effective.
	\end{itemize}
	
	Experiments are performed on a titanium-based alloy under laboratory conditions, which show good agreement with the preliminary calculations. 
	The predicted forces match the predictions of Young for small chip thickness and deviates from it when pure geometric assumptions are less reliable. 
	The predictions of the compared models are similar at small uncut chip thickness and diverge for high feeds and extreme uncut chip geometries.
	The experiments also revealed that the radial cutting force in turning is the most sensitive to modeling assumptions, which is of high importance for the proper development of a dynamical model for turning.
	The uncut curved chip thickness model provided the best prediction for the radial cutting force.
}

\section*{Acknowledgement}
The research reported in this paper has been supported by the project EUROSTARS FORTH E!12998 and by the Ministry for Innovation and Technology (2019-2.1.2-NEMZ-2019-00005).

This research work has been done under the framework of the
project MIRAGED: Posicionamiento estratégico en modelos virtuales y gemelos digitales para una industria 4.0 (CER-20191001), supported by CDTI-acreditación y concesión de ayudas destinadas a centros tecnológicos de excelencia cervera.



\section*{Appendix}
\appendix
\section{Finite element model}\label{app:FE}

A simplified finite element model (FE) is presented to replace complicated and time-consuming numerical simulations. The gradient of the displacement field is found to be applicable to predict the local chip flow directions on the reference plane. The details on the FE model are presented below. 
Description on finite element models are found in many textbooks, for an introductory course, see the book of Steven M. Lepi \cite{Lepi1998book}.

\begin{figure}[!tb]
	\centering
	\includegraphics[scale=1]{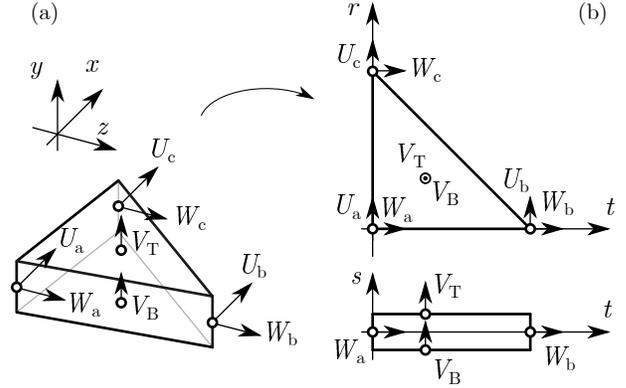}
	\caption{{\bf Triangular finite element used for compression.} (a) Global coordinate system. (b) Local coordinate system.}
	\label{fig:FE_element}
\end{figure}

To speed up numerical calculation, a wedge-shaped finite element is applied with reduced number of degrees of freedom. The element is shown in Fig.~\ref{fig:FE_element}(a-b), which has eight degrees of freedom in total.
The thickness is constant, and the external load is compression only normal to the plane. 
The element is parallel with the ($x,z$) plane, the vertices are denoted by $\rm a$, $\rm b$, $\rm c$, while the bottom and top layers by $\rm B$ and $\rm T$. 
The nodal displacements are therefore $U_{\rm a}$, $W_{\rm a}$, $U_{\rm b}$, $W_{\rm b}$, $U_{\rm c}$, $W_{\rm c}$, $V_{\rm B}$ and  $V_{\rm T}$.

The transformation between the local ($x,y,z$) and natural coordinate system ($r,s,t$) is described by the interpolation functions $N_i$ as
\begin{equation}\label{eq:interp1}
	\begin{split}
		{x}(r,s,t)& =N_1x_{\rm a}+N_2x_{\rm b}+N_3x_{\rm c},\\
		{y}(r,s,t)& =N_4y_{\rm B}+N_5y_{\rm T},\\
		{z}(r,s,t)& =N_1z_{\rm a}+N_2z_{\rm b}+N_3z_{\rm c},
	\end{split}
\end{equation}
where
\begin{equation}
	\begin{split}
		N_1 & =(1-r-t), \quad
		N_2=t, \quad
		N_3=r, \\
		N_4 & =1-s, \quad
		N_5=s
	\end{split}
\end{equation}
are the shape functions. Note, that the shape functions are linear, the derivatives are constants, so the interpolated stresses and strains are constant inside the element.
The same shape functions are used for the interpolation of the displacement fields, i.e.,
\begin{equation}
	\begin{split}
		\tilde{U}(r,s,t)& =N_1U_{\rm a}+N_2U_{\rm b}+N_3U_{\rm c},\\
		\tilde{V}(r,s,t)& =N_4V_{\rm B}+N_5V_{\rm T},\\
		\tilde{W}(r,s,t)& =N_1W_{\rm a}+N_2W_{\rm b}+N_3W_{\rm c}.
	\end{split}
\end{equation}
{\color{black}
	Following the classical approach of finite element modeling, the elemental stiffness matrix can be calculated as
	\begin{equation}
		{\bf K}^{\rm e}=\int_{(V)} {\bf B}^\top{\bf E}\,{\bf B}\,{\rm d}V=\frac{1}{2}{\bf B}^\top{\bf E}\,{\bf B}\, {\det({\bf J})},
	\end{equation}
	where $\bf B$ is the strain-displacement matrix, $\bf E$ is the material matrix (containing the modulus of elasticity) and $\bf J$ is the Jacobian matrix \cite{Lepi1998book}.
	Note that the integration over the volume of the element is expressed explicitly due the linear interpolation functions.}

The equilibrium equation for the entire model is formulated as \begin{equation}
	{\bf K}\,{\bf U}={\bf F},
\end{equation}
where $\bf K$, $\bf U$, and $\bf F$ are the global stiffness matrix, displacement vector, and load vector, which are assembled from the elemental matrices. In our model, we compress the plate normal to its plane (perpendicular to $y$), which is modeled as a constraint (displacement-driven load). 
Therefore, the elements of $\bf U$ contain not only zeros and unknowns, but also the prescribed displacements for $V_k$ ($k=1,\dots,N_{\rm e}$, $N_{\rm e}$ is the number of the elements), e.g., 1 and 0 at the top and bottom coordinates, respectively.  
The inplane deformations $U_{l}$ and $W_l$ are the unknowns ($l=1,\dots,N_{\rm n}$, $N_{\rm n}$ is the number of the nodes), moreover, the nodes located on the cutting edge are also constrained in direction $x$ and $z$. Since the equation is linear, the solution is proportional to the load, however, only the gradient of the solution is used, the value of the compression (and material parameters, such as $E$) is arbitrary. 

The nodal displacements $U_l$ and $W_l$ are used to determine the gradient of the displacement field.
The norm of the displacement (ignoring the compression in direction $y$) is calculated as
\begin{equation}
	{G}_i=\sqrt{{U}_i^2+{W}_i^2}, \quad i={\rm a, b, c}, \\
\end{equation}
and the same interpolation functions are used, i.e.,
\begin{equation}
	\tilde{G}(r,s,t) =N_1G_{\rm a}+N_2G_{\rm b}+N_3G_{\rm c}.
\end{equation}
Using the same procedure, the gradient is calculated from the shape functions and Jacobian matrix as
\begin{equation}
	{\bf grad} \tilde G={\setstretch{2.0}{\bf J}^{-1}\begin{bmatrix}
			\ffrac{\partial N_1}{\partial r} &\ffrac{\partial N_2}{\partial r}  &\ffrac{\partial N_3}{\partial r} \\
			0&0 &0 \\
			\ffrac{\partial N_1}{\partial t}  &\ffrac{\partial N_2}{\partial t} &\ffrac{\partial N_3}{\partial t}
	\end{bmatrix}}\begin{bmatrix}
		{G}_{\rm a}\\{G}_{\rm b}\\{G}_{\rm c}
	\end{bmatrix}.
\end{equation}
Finally, the vector field can be normalized, and only the projected inplane components are left, i.e.,

\begin{equation}
	{\bf g}(x,z)=\begin{bmatrix}
		g_x(x,z)\\0\\g_z(x,z)
	\end{bmatrix}:=\frac{	{\bf grad} \tilde G}{|	{\bf grad} \tilde G|},
\end{equation}
where the component in $y$ is zero due to the simplifications. 
The artificial displacement gradient ${\bf g}(x,z)$ defines the direction of the chip flow on the reference plane.

\small

%
\end{document}